\documentclass[11pt]{article}

\usepackage[draft]{ed}
\usepackage{amsfonts}
\usepackage{amsmath}
\usepackage[linesnumbered,lined,boxed,commentsnumbered]{algorithm2e}
\usepackage{amssymb}

\usepackage{epsfig}
\usepackage{bm}
\usepackage{enumerate}
\numberwithin{equation}{section} 

\usepackage{color}
\usepackage[normalem]{ulem}

\title{Combined field-only boundary integral equations for PEC electromagnetic scattering problem in spherical geometries}
\author{Luiz Faria \thanks{POEMS Laboratory (CNRS/INRIA/ENSTA Paris), Institut Polytechnique de Paris, 91120 Palaiseau, France, e-mail: luiz.maltez-faria@inria.fr.}\and Carlos P\'erez-Arancibia \thanks{Department of Applied Mathematics, University of Twente, Enschede, The Netherlands, e-mail:
c.a.perezarancibia@utwente.nl.} \and Catalin Turc\thanks{  Department of
Mathematical Sciences, New Jersey  Institute of Technology,
Univ. Heights. 323 Dr. M. L. King Jr. Blvd, Newark, NJ 07102, USA, e-mail: catalin.c.turc@njit.edu.}}

\newtheorem{theorem}{Theorem}[section]
\newtheorem{lemma}[theorem]{Lemma}

\newtheorem{remark}[theorem]{Remark}
\newenvironment{proof}{\hspace{0.5cm} {\bf Proof.}}
{$\quad {}_\blacksquare$\vspace{0.3cm}}

\date{}
\newcommand{\triple}[1]{{\left\vert\kern-0.25ex\left\vert\kern-0.25ex\left\vert #1 
    \right\vert\kern-0.25ex\right\vert\kern-0.25ex\right\vert}}

\begin{document}
\maketitle
\begin{abstract}
  We analyze the well posedness of certain field-only boundary integral
  equations (BIE) for frequency domain electromagnetic scattering from perfectly
  conducting spheres. Starting from the observations that (1) the three
  components of the scattered electric field $\mathbf{E}^s(\mathbf{x})$ and (2)
  scalar quantity $\mathbf{E}^s(\mathbf{x})\cdot\mathbf{x}$  are radiative
  solutions of the Helmholtz equation, novel boundary integral equation
  formulations of electromagnetic scattering from perfectly conducting obstacles
  can be derived using Green's identities applied to the aforementioned
  quantities and the boundary conditions on the surface of the scatterer. The
  unknowns of these formulations are the normal derivatives of the three
  components of the scattered electric field and the normal component of the
  scattered electric field on the surface of the scatterer, and thus these
  formulations are referred to as field-only BIE. In this paper we use the
  Combined Field methodology of Burton and Miller within the field-only BIE
  approach and we derive new boundary integral formulations that feature only
  Helmholtz boundary integral operators, which we subsequently show to be well
  posed for all positive frequencies in the case of spherical scatterers.
  Relying on the spectral properties of Helmholtz boundary integral operators in
  spherical geometries, we show that the combined field-only boundary integral
  operators are diagonalizable in the case of spherical geometries and their
  eigenvalues are non zero for all frequencies. Furthermore, we show that for
  spherical geometries one of the field-only integral formulations considered in
  this paper exhibits eigenvalues clustering at one---a property similar to
  second kind integral equations.
   
\end{abstract}

\section{Introduction}
\label{intro}

The Maxwell scattering problem in the frequency domain with homogeneous material parameters and perfectly electrically conducting (PEC) boundary conditions can be formulated in an equivalent manner through Boundary Integral Equations~\cite{HsiaoKleinman,Nedelec}.  In the case when the surface of the scatterer is a closed surface in three dimensions, the most widely used strategy to derive BIE formulations relies on Stratton-Chu~\cite{HsiaoKleinman} representation formulas that leads to the classical Magnetic Field and Electric Field Integral Equations (MFIE and EFIE respectively) whose unknown is  the magnetic current (a tangential vector field) on the surface of the scatterer. These two formulations can be combined in a formulation (referred to as the Combined Field Integral Equation) that is well posed for all frequencies~\cite{HarringtonMautz}, and as such it is the robust formulation of choice for most integral equation based numerical methods for the scattering simulations. Recent contributions~\cite{AlougesLevadoux,br-turc,turc3} extended the combined field technology to derive regularized combined field integral equations of the second kind for Maxwell scattering problems. A different approach to derive robust boundary integral equation formulations for Maxwell scattering problems is based on the generalized Debye source representation of electromagnetic fields~\cite{epstein}. The Debye sources can be expressed in terms of the surface divergence of the electric and magnetic currents on the surface of the scatterer, and these two scalar surface functional densities are the unknowns of the Maxwell BIEs.

Following the contributions~\cite{Klaseboer,Yuffa1,Yuffa2}, we study in this
paper a different type of Maxwell BIE formulations whose unknowns, rather than
being surface currents or their charges, are certain boundary values of the
scattered electric fields on the surface of the scatterer. For the latter
reason, these formulations are referred to as field-only. The main advantage of
these formulations is their exclusive reliance on Helmholtz Boundary Integral
Operators, which are easier to discretize than the electromagnetic BIOs in both
the Galerkin and Nystr\"om framework. The trade-off here is the fact that
field-only formulations require three or four surface unknown densities rather
than two surface functional unknowns in the classical Maxwell BIEs.
The
field-only Maxwell BIEs use Helmholtz Green's identities for the three scalar
components of the scattered electric field  $\mathbf{E}^s(\mathbf{x})$ and the
scalar quantity $\mathbf{E}^s(\mathbf{x})\cdot\mathbf{x}$, the latter also being
a radiative solution of the Helmholtz equation in the exterior of a closed
scatterer. Using the PEC boundary conditions, the boundary values of the
scattered electric field on the scatterer can be expressed in terms of the
normal derivatives of the scattered electric field $\partial_n\mathbf{E}^s$ and
its normal component $\mathbf{E}^s\cdot\mathbf{n}$, where $\mathbf{n}$ is the
outward unit normal on the surface of the PEC scatterer. Therefore, the
application of the boundary Dirichlet and Neumann traces on the aforementioned
Green's identities leads to a system of four Helmholtz BIEs in terms of the four
field-only unknown quantities $\{\partial_n\mathbf{E}^s,
\mathbf{E}^s\cdot\mathbf{n}\}$. In addition, it is possible to use yet again the
PEC boundary conditions to express $\mathbf{E}^s\cdot\mathbf{n}$ in terms of the
three components of the vector field-only quantity $\partial_n\mathbf{E}^s$ and
the mean-curvature on the scatterer, and thus to derive a system of three
Helmholtz BIE for the solution of Maxwell scattering problems. Unlike the
approach in~\cite{Klaseboer,Yuffa1,Yuffa2}, we pursue the classical combined
field approach of Burton-Miller and
Brackhage-Werner~\cite{BrackhageWerner,BurtonMiller} to the Dirichlet and
Neumann traces applied to Helmholtz Green's identities recounted above. Our
approach leads to systems of BIEs that can be expressed in operator form as
rank-one perturbation of CFIE BIOs that are invertible for all positive
frequencies. Therefore, the invertibility of these BIE systems can be reduced
via the Woodbury formula~\cite{Woodbury} to that of a scalar (albeit
non-standard) BIE. This paper studies the spectral properties of the latter
scalar BIE in the case of spherical geometries, for which these BIEs are
invertible for all positive frequencies. The unique solvability of these
formulations for general scatterers remains an open question. In the case of
spherical geometries, the eigenvalues of the BIE formulation that uses the
unknowns $\partial_n\mathbf{E}^s$ accumulate at infinity, while, remarkably, the
eigenvalues of the augmented BIE formulation that uses the unkonwns
$\{\partial_n\mathbf{E}^s, \mathbf{E}^s\cdot\mathbf{n}\}$ accumulate at $1$.

The paper is organized as follows: in Section~\ref{sec1} we present the general combined field approach for the derivation of fields-only BIE formulations of Maxwell scattering problems; in Section~\ref{aug} we present an augmented system of BIEs for the solution of Maxwell scattering problems, and finally in Section~\ref{sph} we analyze the spectral properties of these formulations in the case of spherical geometries.


\section{Combined field-only BIE}\label{sec1}

We consider the problem of evaluating the scattered electromagnetic
field $(\mathbf E^{s},\mathbf H^{s})$ that results as an incident
field $(\mathbf E^{i},\mathbf H^{i})$ impinges upon the boundary
$\Gamma$ of a perfectly conducting scatterer $\Omega$. The total
field $(\mathbf E,\mathbf H) = (\mathbf E^{s}+\mathbf E^{i},\mathbf H^{s}+\mathbf H^{i})$ satisfies the time harmonic Maxwell equations with wavenumber $k>0$
\begin{equation}
  \label{eq:Maxwell}
  {\rm curl}\ {\mathbf E}-ik {\mathbf H}=\mathbf{0},\qquad {\rm curl}\ {\mathbf H}+ik\mathbf{E}=\mathbf{0} \qquad \rm{in}\ \Omega^+:=\mathbf{R}^{3}\setminus \Omega
\end{equation}
together with the perfect-conductor boundary conditions
\begin{equation}
\label{eq:bc}
\mathbf{n}\times\mathbf{E}=\mathbf{0}\qquad \rm{on}\ \Gamma.
\end{equation}
Here and in what follows $\mathbf{n}$ denotes the unit normal field on $\Gamma$ pointing into $\Omega^+$. Also, for each $\mathbf{x}\in\Gamma$, the unit vector $\mathbf{n}(\mathbf{x})$ will denote a column vector; furthermore, the electric and magnetic fields are column vectors. The incident field $(\mathbf E^{i},\mathbf H^{i})$ satisfies the Maxwell equations. In addition, in order to ensure well posedness of the Maxwell equations in the exterior domain $\Omega^+$, we require that the scattered field $(\mathbf E^{s},\mathbf H^{s})$ satisfies the well known Silver-M\"uller radiation conditions at infinity~\cite{Nedelec}. The components of the scattered electric field, which we denote by $E^s_j, j=1,2,3$, satisfy in turn the Helmholtz equation with the same wavenumber $k$, that is
\[
\Delta u+k^2 u=0\qquad \rm {in}\ \Omega^+
\]
and also the Sommerfeld radiation conditions at infinity~\cite{Nedelec}. Consequently, the components of the scattered electric fields can be represented using Green's identities in the domain $\Omega^+$
\begin{equation}\label{eq:id_G0}
  E^s_j(\mathbf{x})=DL[(E^s_j)|_\Gamma](\mathbf{x})-SL[(\partial_nE^s_j)|_\Gamma](\mathbf{x}),\ \mathbf{x}\in\Omega^+,\ j=1,2,3\\
\end{equation}
in terms of their Dirichlet and Neumann boundary values on $\Gamma$ and the single and double layer potentials associated with the radiative Helmholtz Green's function. Specifically, for functional densities $\varphi$ and $\psi$ defined on $\Gamma$, the Helmholtz layer potentials are defined in the following manner
\[
DL[\psi](\mathbf{x}):=\int_\Gamma \frac{\partial G_k(\mathbf{x}-\mathbf{y})}{\partial \mathbf{n}(\mathbf{y})}\ \psi(\mathbf{y}) ds(\mathbf{y})\quad SL[\varphi](\mathbf{x}):=\int_\Gamma G_k(\mathbf{x}-\mathbf{y}) \varphi(\mathbf{y})ds(\mathbf{y}),\quad \mathbf{x}\in\mathbb{R}^3\setminus\Gamma
\]
in terms of the Helmholtz Green's function corresponding to the wavenumber $k$
\[
G_k(\mathbf{x}):=\frac{e^{ik|\mathbf{x}|}}{4\pi|\mathbf{x}|}.
\]
We will derive boundary integral equation (BIE) formulations for the Maxwell scattering problem by applying Dirichlet and Neumann boundary conditions on $\Gamma$ to the Green's identities~\eqref{eq:id_G0}. The unknowns of these BIE are the Neumann boundary data $(\partial_nE^s_j)|_\Gamma$ for $j=1,2,3$. To this end, we start by expressing the Dirichlet boundary data  $(E^s_j)|_\Gamma$ for $j=1,2,3$ in terms of the Neumann boundary data. First, we use the PEC boundary condition to express $(E^s_j)|_\Gamma$ for $j=1,2,3$ in terms of the normal component of the electric field on the boundary $(\mathbf{E}^s\cdot\mathbf{n})|_\Gamma$ and the incoming electric field. Indeed, we use the last two scalar components PEC boundary conditions and we get
\begin{eqnarray*}
  n_2E^s_1-n_1E^s_2&=&n_1E^i_2-n_2E^i_1\\
  n_3E^s_1-n_1E^s_3&=&n_1E^i_3-n_3E^i_1\\
  n_1E^s_1+n_2E^s_2+n_3E^s_3&=&\mathbf{E}^s\cdot\mathbf{n}
\end{eqnarray*}
where $\mathbf{n}=(n_1,n_2,n_3)$. Multiplying both sides of the first equation above by $n_2$, the second equation by $n_3$, and the third equation by $n_1$ and adding the ensuing results, we obtain immediately an expression of the first component of the Dirichlet boundary value of scattered electric field on $\Gamma$ in terms of the normal component of the scattered field on the boundary
\begin{equation*}
  E^s_1=n_1\mathbf{E}^s\cdot\mathbf{n}+(\mathbf{E}^i_t)_1
\end{equation*}
where $\mathbf{E}^i_t:=\mathbf{n}\times(\mathbf{n}\times\mathbf{E}^i)$. Similarly, we obtain
\begin{equation}\label{eq:first_id}
  E^s_j=n_jE^s_n+(\mathbf{E}^i_t)_j,\ j=1,2,3.
\end{equation}
where we have denoted $E^s_n:=(\mathbf{E}^s\cdot\mathbf{n})|_\Gamma$. Second, we connect the quantities $E^s_n$ and $(\partial_n \mathbf{E}^s)|_\Gamma$ via certain differential geometry identities. To this end, using the PEC boundary conditions~\eqref{eq:bc} we get that
\begin{equation}\label{eq:totalE}
\mathbf{E}|_\Gamma=(\mathbf{E}\cdot\mathbf{n})|_\Gamma\ \mathbf{n}.
\end{equation}
We derive next an on surface relation between the quantities $\partial_n \mathbf{E}$ and $\mathbf{E}\cdot\mathbf{n}$. To this end, we make use of certain differential geometry identities following the reference~\cite{Nedelec}. First, we define the tubular neighborhood $\Gamma_\varepsilon$ of $\Gamma$ as the set of all points whose distance to the surface $\Gamma$ is less than $\varepsilon$. For $C^2$ surfaces $\Gamma$ and small enough $\varepsilon$, each point $\mathbf{y}\in\Gamma_\varepsilon$ has a unique projection onto $\Gamma$ which we denote by $\mathcal{P}(\mathbf{y})$. In this manner, any point $\mathbf{y}\in\Gamma_\varepsilon$ can be expressed in the form $\mathbf{y}=\mathcal{P}(\mathbf{y})+s \mathbf{n}(\mathcal{P}(\mathbf{y})),\ -\varepsilon\leq s\leq\varepsilon$. For a vector field defined on $\Gamma_\varepsilon$, we use the following decomposition of its divergence~\cite{Nedelec} along parallel surfaces
\begin{equation}\label{eq:div}
\nabla\cdot\mathbf{v}=\nabla_{\Gamma_s}\cdot \mathbf{v}_{\Gamma_s}+2\mathcal{H}(\mathbf{v}\cdot\mathbf{n})+\frac{\partial}{\partial s}(\mathbf{v}\cdot\mathbf{n}),\qquad \mathbf{v}_{\Gamma_s}:=\mathbf{n}\times(\mathbf{v}\times\mathbf{n})
\end{equation}
where $\mathcal{H}$ is the mean curvature on $\Gamma$ and $\nabla_{\Gamma_s}$ denotes the tangential gradient on the parallel surface $\Gamma_s$. Using the fact that the electric field is divergence free and that the curvature operator $\nabla \mathbf{n}$ is a tangential operator, and taking into account equation~\eqref{eq:totalE}, we obtain from the decompostion~\eqref{eq:div} the following relation on $\Gamma$ 
\begin{equation}\label{eq:curvature}
  \partial_n\mathbf{E}\cdot \mathbf{n}=-2\mathcal{H}\ \mathbf{E}\cdot \mathbf{n}.
\end{equation}
Thus, the normal component of the scattered electric field on $\Gamma$ can be connected with the normal derivative of the scattered electric field on $\Gamma$ through the relation
\begin{equation}\label{eq:curvC}
  E_n^s=-\frac{1}{2}\mathcal{H}^{-1}\partial_n \mathbf{E}^s\cdot \mathbf{n}-\frac{1}{2}\mathcal{H}^{-1}\partial_n \mathbf{E}^i\cdot \mathbf{n}-E^i_n.
\end{equation}
Combining equations~\eqref{eq:first_id} and~\eqref{eq:curvC} we obtain
\begin{equation}\label{eq:First_id}
  (\mathbf{E}^s)|_\Gamma=- \frac{1}{2}\mathcal{H}^{-1}\mathbf{n}\mathbf{n}^\top (\partial_n \mathbf{E}^s)|_\Gamma- \frac{1}{2}\mathcal{H}^{-1}\mathbf{n}\mathbf{n}^\top (\partial_n \mathbf{E}^i)|_\Gamma-(\mathbf{E}^i)|_\Gamma,
\end{equation}
giving the desired relation between the scattered electric field $\mathbf{E}^s$ 
and its normal derivative $\partial_n\mathbf{E}^s$ on $\Gamma$.
We are now in the position to derive a system of BIE for the solution of the Maxwell scattering problems whose unknown is the vector quantity $(\partial_n\mathbf{E}^s)|\Gamma$. Applying both Dirichlet and Neumann traces to the Green's identities~\eqref{eq:id_G0} we derive the following vector integral relations
\begin{equation}\label{eq:G1}
  \frac{1}{2}\mathbf{E}^s - K[\mathbf{E}^s] + S[\partial_n \mathbf{E}^s]=0\quad {\rm on}\ \Gamma
\end{equation}
and respectively
\begin{equation}\label{eq:G2}
  \frac{1}{2}\partial_n\mathbf{E}^s + K^\top[\partial_n\mathbf{E}^s] -N[\mathbf{E}^s]=0\quad {\rm on}\ \Gamma.
\end{equation}
In equations~\eqref{eq:G1} and~\eqref{eq:G2} the boundary integral operators $S, K, K^\top, N$ are the four boundary integral operators associated with the Calder\'on's calculus for the Helmholtz equation. For two scalar functional densities $\varphi$ and $\psi$ defined on $\Gamma$, the four Helmholtz boundary integral operators are defined as for $\mathbf{x}\in\Gamma$
\begin{eqnarray*}
  S[\varphi](\mathbf{x}):=\int_\Gamma G_k(\mathbf{x}-\mathbf{y}) \varphi(\mathbf{y})ds(\mathbf{y}) && K[\psi](\mathbf{x}):=\int_\Gamma \frac{\partial G_k(\mathbf{x}-\mathbf{y})}{\partial \mathbf{n}(\mathbf{y})}\ \psi(\mathbf{y}) ds(\mathbf{y})\\
  K^\top[\varphi](\mathbf{x}):=\int_\Gamma \frac{\partial G_k(\mathbf{x}-\mathbf{y})}{\partial \mathbf{n}(\mathbf{x})}\ \varphi(\mathbf{y}) ds(\mathbf{y}) && N[\psi](\mathbf{x}):=f.p. \int_\Gamma \frac{\partial^2 G_k(\mathbf{x}-\mathbf{y})}{\partial \mathbf{n}(\mathbf{x})\partial \mathbf{n}(\mathbf{y})}\ \psi(\mathbf{y}) ds(\mathbf{y}).
\end{eqnarray*}
We remark that in equations~\eqref{eq:G1} and~\eqref{eq:G2} the application of the boundary integral operators to vector functional densities is performed componentwise. We combine classically equations~\eqref{eq:G1} and~\eqref{eq:G2} in the Burton-Miller and Brackhage-Werner manner~\cite{BrackhageWerner,BurtonMiller}. Choosing a coupling parameter $\eta>0$, we subtract from equation~\eqref{eq:G2} the equation~\eqref{eq:G1} multiplied by the factor $i\eta$ and we derive the integral identity
\begin{equation}\label{eq:id0}
\mathcal{D}(\eta)[\partial_n \mathbf{E}^s]-i\eta\ \mathcal{N}(\eta^{-1})[\mathbf{E}^s]=0\quad {\rm on}\ \Gamma
\end{equation}
where we have used the following notations
\begin{equation}\label{eq:ops}
  \mathcal{D}(\eta):=\frac{1}{2}I+K^\top-i\eta\ S\qquad \mathcal{N}(\eta):=\frac{1}{2}I-K-i\eta\ N,\ \eta\neq 0.
  \end{equation}
Finally, using relation~\eqref{eq:First_id} to replace $\mathbf{E}^s$ in favor
of $\partial_n \mathbf{E}^s$ in equations~\eqref{eq:id0}, we derive the following system of BIE on $\Gamma$
\begin{equation}\label{eq:3x3}
  \boldsymbol{\mathcal{A}}_1[\partial_n\mathbf{E}^s]=-i\eta\mathcal{N}(\eta^{-1})\left[E^i+\frac{1}{2}\mathcal{H}^{-1}\mathbf{n}\mathbf{n}^\top\partial_n\mathbf{E}^i\right]
\end{equation}
where the boundary integral operator $\boldsymbol{\mathcal{A}}_1$ is defined as 
\begin{equation}
  \boldsymbol{\mathcal{A}}_1:=\mathcal{D}(\eta)\mathbb{I}_3+i\eta\mathcal{N}(\eta^{-1})\left[\frac{1}{2}\mathcal{H}^{-1}\mathbf{n}\mathbf{n}^\top\right].
\end{equation}
Clearly, taking into account the mapping properties of the four Helmholtz BIO, we can see that the matrix operator $\boldsymbol{\mathcal{A}}_1$ has the following mapping property $\boldsymbol{\mathcal{A}}_1:(H^1(\Gamma))^3\to (L^2(\Gamma))^3$, where $H^s(\Gamma)$ denotes the Sobolev space of order $s$ on $\Gamma$~\cite{Nedelec}. Given that the operator $\mathcal{D}(\eta):H^s(\Gamma)\to H^s(\Gamma)$ is an invertible Brackage-Werner type combined field operator~\cite{Nedelec} for all Sobolev index $s$, the invertibility of the integral operator $\boldsymbol{\mathcal{A}}_1$ is equivalent to that of the integral operator $\boldsymbol{\mathcal{B}}_1:(H^1(\Gamma))^3\to (L^2(\Gamma))^3$ defined as 
\[
\boldsymbol{\mathcal{B}}_1=\mathbb{I}_3+\mathcal{S}(\eta)\left[\frac{1}{2}\mathcal{H}^{-1}\mathbf{n}\mathbf{n}^\top\right],\quad \mathcal{S}(\eta):=i\eta\mathcal{D}^{-1}(\eta)\mathcal{N}(\eta^{-1}).
\]
Using the fact that the matrix operator $\mathbf{n}\mathbf{n}^\top$ is a rank-one operator, we can reduce the question of the invertibility of the operator $\boldsymbol{\mathcal{B}}_1$ to that of a scalar operator. Indeed, solving the vector equation $\boldsymbol{\mathcal{B}}_1 \boldsymbol{\psi}=\mathbf{f}$ can be written in the form
\begin{equation}\label{eq:vec}
  \boldsymbol{\psi}+\mathcal{S}(\eta)\left[\frac{1}{2}\mathcal{H}^{-1}\mathbf{n}\left(\mathbf{n}\cdot\boldsymbol{\psi}\right)\right]=\mathbf{f}.
\end{equation}
Multiplying both sides of equation~\eqref{eq:vec} by $\mathbf{n}^\top$ we obtain
the following scalar BIE
\begin{equation}\label{eq:vec1}
  \psi_n+\mathbf{n}\cdot\mathcal{S}(\eta)\left[\frac{1}{2}\mathcal{H}^{-1}\mathbf{n}\ \psi_n\right]=\mathbf{f}\cdot\mathbf{n},\quad \psi_n:=\mathbf{n}\cdot\boldsymbol{\psi}\in H^1(\Gamma).
\end{equation}

 Note that~\eqref{eq:vec1} suggests a way of solving~\eqref{eq:3x3} by
considering only decoupled scalar problems. In particular, having found $\psi_n
= \boldsymbol{\psi} \cdot \mathbf{n}$ by means of~\eqref{eq:vec1}, the
tangential components of the Neumann trace $\partial_n \mathbf{E}^s$ are simply
given by
\begin{equation}\label{eq:tan-comp-vec1}
  \mathbf{t}_i \cdot \boldsymbol{\psi} = \mathbf{t}_i \cdot \left( \mathbf{f} - \mathcal{S}(\eta)\left[\frac{1}{2}\mathcal{H}^{-1}\mathbf{n}\
  \psi_n\right] \right),\quad i = 1,2,
\end{equation}
where $\mathbf{t}_1(\boldsymbol{x})$ and $\mathbf{t}_2(\boldsymbol{x})$ form a
basis for the tangent plane at $\boldsymbol{x} \in \Gamma$. Thus, despite having
$3$ surface unknowns in $\boldsymbol{\psi}$, the Neumann derivative of the
electric field can be found by solving three scalar problems. 

We will investigate in Section~\ref{sph} the unique solvability of the
BIE~\eqref{eq:vec1} for spherical geometries. 

\begin{remark}[Evaluation of $\mathcal{S}$]
The evaluation of $\mathcal{S} : H^1(\Gamma) \to L^2(\Gamma)$ may be performed
iteratively by solving a sound-soft Helmholtz scattering problem on $\Gamma$
using the combined field formulation. More precisely, letting $\sigma \in
H^1(\Gamma)$ be a given surface density, we have that $g =
\mathcal{S}(\eta)[\sigma]$ solves 
\begin{align}
  \mathcal{D}(\eta)[g] = i\eta \mathcal{N}(\eta^{-1})[\sigma].
\end{align}
Since $\mathcal{D}$ is a second-kind integral operator, we expect the above
system to be well-suited to iterative solvers such as GMRES. This means, in
particular, that the inverse of $\mathcal{D}$ does not need to be explicitly
constructed if one wishes to solve~\eqref{eq:vec1}.
\end{remark}

\section{An augmented system of scalar BIE for Maxwell's scattering problems}\label{aug}

We derive in this section a system of $4\times 4 $ scalar BIE for Maxwell scattering problems. The main idea is to augment the use of Green's identity for the three components of the scattered electric field $E^s_j,j=1,2,3$ with the Green's identities  for the scalar quantity $\mathbf{E}^s(\mathbf{x})\cdot\mathbf{x}$. A simple calculation shows that $\mathbf{E}^s(\mathbf{x})\cdot\mathbf{x}$ satisfies the Helmholtz equation in the domain $\Omega^+$. We prove in what follows that $\mathbf{E}^s(\mathbf{x})\cdot \mathbf{x}$ is actually a {\em radiating} solution of the Helmholtz equation in $\Omega^+$. To that end, we begin with the representation of the radiating field $\mathbf{E}^s(\mathbf{x})$ for $|\mathbf{x}|\geq R$ where $R$ is  large enough such that $\overline{\Omega}\subset B_R(0)$ [~\cite{Nedelec} equations (5.3.2) and (5.3.4)]:
\begin{eqnarray}\label{eq:radiationE}
  \mathbf{E}^s(\mathbf{x})&=&\sum_{\ell=1}^{\infty}\sum_{m=-\ell}^\ell a_{\ell}^m h_\ell^{(1)}(kr)T_\ell^m(\hat{\mathbf{x}})\nonumber\\
  &+&i\sum_{\ell=1}^\infty\sum_{m=-\ell}^\ell \frac{b_\ell^m}{2\ell+1}\left[(\ell+1)h_{\ell-1}^{(1)}(kr)I_{\ell-1}^m(\hat{\mathbf{x}})+\ell h_{\ell+1}^{(1)}(kr)N_{\ell+1}^m(\hat{\mathbf{x}})\right],\ r:=|\mathbf{x}|,\ \hat{\mathbf{x}}:=\mathbf{x}/r\nonumber\\
  \end{eqnarray}
where [\cite{Nedelec} equations (2.4.172), (2.4.173), and (2.4.174)]
\begin{eqnarray}\label{eq:spherical_harm}
  I_\ell^m(\hat{\mathbf{x}})&:=&\nabla_S Y_{\ell+1}^m(\hat{\mathbf{x}})+(\ell+1)Y_{\ell+1}^m(\hat{\mathbf{x}})\hat{\mathbf{x}}\\
  T_\ell^m(\hat{\mathbf{x}})&:=&\nabla_S Y_{\ell}^m(\hat{\mathbf{x}})\times \hat{\mathbf{x}}\\
  N_\ell^m(\hat{\mathbf{x}})&:=&-\nabla_S Y_{\ell-1}^m(\hat{\mathbf{x}})+\ell Y_{\ell-1}^m(\hat{\mathbf{x}})\hat{\mathbf{x}},
\end{eqnarray}
and $h_\ell^{(1)}$ denote the spherical Hankel functions and $Y_\ell^m$ denote the classical spherical harmonics of degree $\ell$ and order $m$. More precisely, here in what follows we will use the following definition of spherical harmonics: given a unit vector $\hat{x}$ on the unit sphere $\mathbb{S}^2$ in the form $\hat{x}=(\cos{\varphi}\sin{\theta},\sin{\varphi}\sin{\theta},\cos{\theta})$, the spherical harmonics of degree $\ell\geq 0$ and order $m, -\ell\leq m\leq\ell$ are defined as
  \[
  Y_\ell^m(\hat{x})=\gamma_\ell^m e^{im\varphi}\ \mathbb{P}_\ell^{|m|}(\cos{\theta}),\qquad \gamma_\ell^m:=\left[\frac{\left(\ell+\frac{1}{2}\right)}{2\pi}\frac{(\ell-|m|)!}{(\ell+|m|)!}\right]^{\frac{1}{2}}
  \]
  where
  \[
  \mathbb{P}_\ell^m(\cos{\theta})=(-1)^m(\sin{\theta})^m\left(\frac{d}{dx}\right)^m \mathbb{P}_\ell(\cos{\theta}),\ 0\leq m\leq \ell.
  \]
  Here $\mathbb{P}_\ell$ are the standard Legendre polynomials. With this definition, the family of spherical harmonics $\{Y_\ell^m\}_{0\leq \ell, -\ell\leq m\leq\ell}$ constitutes an orthonormal basis of $L^2(\mathbb{S}^2)$. It follows then that for all $|\mathbf{x}|\geq R$ we have that
\begin{equation*}\label{eq:rad_Edotx}
  \mathbf{E}^s(\mathbf{x})\cdot\mathbf{x}=i\sum_{\ell=1}^\infty\sum_{m=-\ell}^\ell \frac{b_\ell^m}{2\ell+1}\ell(\ell+1)\left[h_{\ell-1}^{(1)}(kr)+h_{\ell+1}^{(1)}(kr)\right]rY_{\ell}^m(\hat{\mathbf{x}}),
\end{equation*}
if we take into account the fact that $\nabla_S Y_p^q(\hat{\mathbf{x}})\cdot \hat{\mathbf{x}}=0$. Given that [~\cite{Nedelec} (2.6.11)]
\[
h_{\ell-1}^{(1)}(kr)+h_{\ell+1}^{(1)}(kr)=\frac{2\ell+1}{kr}h_\ell^{(1)}(kr)
\]
we obtain
\begin{equation}\label{eq:rad_Edotx_1}
  \mathbf{E}^s(\mathbf{x})\cdot\mathbf{x}=\frac{i}{k}\sum_{\ell=1}^\infty\sum_{m=-\ell}^\ell b_\ell^m\ell(\ell+1)h_{\ell}^{(1)}(kr)Y_{\ell}^m(\hat{\mathbf{x}}),\ |\mathbf{x}|\geq R
\end{equation}
which confirms that $\mathbf{E}^s(\mathbf{x})\cdot\mathbf{x}$ is a radiating solution of the Helmholtz equation in $\Omega^+$. Therefore, we can apply Green's identities to the quantities $\mathbf{E}^s(\mathbf{x})\cdot\mathbf{x}$ in the domain $\Omega^+$, which are expressed in the form
\begin{equation}\label{eq:id_G}
  \mathbf{E}^s(\mathbf{x})\cdot\mathbf{x}=DL[(\mathbf{E}^s(\mathbf{y})\cdot\mathbf{y})|_\Gamma](\mathbf{x})-SL[\partial_n(\mathbf{E}^s(\mathbf{y})\cdot\mathbf{y})|_\Gamma](\mathbf{x}),\quad \mathbf{x}\in\Omega^+.
\end{equation}
We will use the Green's identities~\eqref{eq:id_G0} and~\eqref{eq:id_G} together with the relations~\eqref{eq:first_id} to derive a system of BIE whose unknowns are $\partial_n\mathbf{E}^s$ and $E^s_n$ on $\Gamma$. The missing ingredient is a relation that expresses the quantity $\partial_n(\mathbf{E}^s\cdot\mathbf{x})|_\Gamma$ in terms of those four scalar unknowns. Such a link is readily obtained via the calculus identity
\begin{equation}\label{eq:second_id}
  \partial_n(\mathbf{E}^s\cdot\mathbf{x})|_\Gamma=\sum_{j=1}^3(x_j\partial_nE^s_j)|_\Gamma+E^s_n.
\end{equation}
First we express the boundary values of the scattered electric field on $\Gamma$ via the identities~\eqref{eq:first_id} in equation~\eqref{eq:G1}, then we apply exterior Dirichlet traces to the Green identities~\eqref{eq:id_G} above and take into account the relations~\eqref{eq:second_id} to arrive at the following system of boundary integral equations
\begin{eqnarray}\label{eq:system_4}
  \frac{1}{2}n_jE^s_n-K\left[n_jE^s_n\right]+S\left[\partial_nE^s_j\right]&=&-\frac{1}{2}(\mathbf{E}^i_t)_j+K\left[(\mathbf{E}^i_t)_j\right],\ j=1,2,3\nonumber\\
  \frac{1}{2}\mathbf{x}\cdot\mathbf{n}E^s_n-K\left[\mathbf{x}\cdot\mathbf{n}E^s_n\right]+S\left[(\sum_{j=1}^3x_j\partial_n E^s_j)+E^s_n\right]&=&-\frac{1}{2}\mathbf{x}\cdot\mathbf{E}^i_t+K\left[\mathbf{x}\cdot\mathbf{E}^i_t\right].
\end{eqnarray}
We repeat the procedure to the Neumann traces equations~\eqref{eq:G1} and we further apply exterior Neumann traces to the Green identities~\eqref{eq:id_G} taking again into account the relations~\eqref{eq:second_id}, and we obtain the following system of boundary integral equations
\begin{eqnarray}\label{eq:system_4N}
  \frac{1}{2}\partial_nE^s_j+K^\top\left[\partial_n E^s_j\right]-N\left[n_j E^s_n\right]&=&N\left[(\mathbf{E}^i_t)_j\right],\ j=1,2,3\nonumber\\
  \frac{1}{2}E_n^s+\frac{1}{2}\sum_{j=1}^3x_j\partial_n E^s_j+K^\top\left[E_n^s+\sum_{j=1}^3x_j\partial_n E^s_j\right]
  -N\left[\mathbf{x}\cdot\mathbf{n}E_n^s\right]&=&N\left[\mathbf{x}\cdot\mathbf{E}^i_t\right].
\end{eqnarray}
We follow the spirit of combined field formulations and we subtract from equations~\eqref{eq:system_4N} those obtained after multiplication by the coupling factor $i\eta,\ \eta>0$ of the equations~\eqref{eq:system_4}. The resulting combined system can be written making use of the following matrix multipliers
\begin{equation*}
  \boldsymbol{\mathcal{M}}_1:=\begin{bmatrix}I & 0 & 0 & 0 \\ 0 & I & 0 & 0\\ 0 & 0 & I & 0\\ x_1 & x_2 & x_3 & I\end{bmatrix}\qquad \boldsymbol{\mathcal{M}}_2:=\begin{bmatrix}0 & 0 & 0 & n_1 \\ 0 & 0 & 0 & n_2\\ 0 & 0 & 0 & n_3\\ 0 & 0 & 0 & 0\end{bmatrix}
  \end{equation*}
in the compact form
\begin{equation}\label{eq:system_m_f}
  \boldsymbol{\mathcal{A}}\begin{bmatrix}\partial_n\mathbf{E}^s\\ E^s_n\end{bmatrix}=i\eta\mathcal{N}(\eta^{-1})\boldsymbol{\mathcal{M}}_1\begin{bmatrix}\mathbf{E}^i_t\\ 0\end{bmatrix}
  \end{equation}
\begin{equation}\label{eq:matrix_form}
  \boldsymbol{\mathcal{A}}:=\mathcal{D}(\eta)\boldsymbol{\mathcal{M}}_1-i\eta\ \mathcal{N}(\eta^{-1})\boldsymbol{\mathcal{M}}_1\boldsymbol{\mathcal{M}}_2
\end{equation}
where the operators $\mathcal{D}$ and $\mathcal{N}$ were defined in equations~\eqref{eq:ops}. Given that the multiplier operator $\boldsymbol{\mathcal{M}}_1:(H^s(\Gamma))^4\to (H^s(\Gamma))^4$ is invertible, the invertibility of the integral operator $\boldsymbol{\mathcal{A}}:(H^1(\Gamma))^4\to (L^2(\Gamma))^4$ is equivalent to that of the operator
\[
\boldsymbol{\mathcal{A}}\boldsymbol{\mathcal{M}}_1^{-1}=\mathcal{D}(\eta)\mathbb{I}_4-i\eta\ \mathcal{N}(\eta^{-1})\boldsymbol{\mathcal{C}},\qquad \boldsymbol{\mathcal{C}}:=\boldsymbol{\mathcal{M}}_1\boldsymbol{\mathcal{M}}_2\boldsymbol{\mathcal{M}}_1^{-1}.
\]
Now, both Helmholtz CFIE operators $\mathcal{D}(\eta):H^s(\Gamma)\to H^s(\Gamma)$ and $\mathcal{N}(\eta^{-1}):H^{s+1}(\Gamma)\to H^s(\Gamma)$ are invertible with continuous inverses for all Sobolev indices $s$, so the invertibility of the operator $\boldsymbol{\mathcal{A}}\boldsymbol{\mathcal{M}}_1^{-1}$ is equivalent in turn to that of the operator $\boldsymbol{\mathcal{B}}:(H^1(\Gamma))^4\to (L^2(\Gamma))^4$ defined as
\begin{equation}\label{eq:inv_f}
  \boldsymbol{\mathcal{B}}:=\mathbb{I}_4-\mathcal{S}(\eta)\boldsymbol{\mathcal{C}}.
\end{equation}
We note that the matrix operator $\boldsymbol{\mathcal{C}}$ is a rank one multiplier
\[
\boldsymbol{\mathcal{C}}=\mathbf{u}\mathbf{v}^\top,\ \mathbf{u}:=\begin{bmatrix}n_1\\n_2\\n_3\\\mathbf{x}\cdot\mathbf{n}\end{bmatrix}\, \mathbf{v}:=\begin{bmatrix}-x_1 \\ -x_2 \\ -x_3 \\ 1\end{bmatrix}.
  \]
  Clearly, both vectors $\mathbf{u}$ and $\mathbf{v}$ depend on $\mathbf{x}\in\Gamma$, but we omit their dependence as we do not believe there is a risk of confusion. As before, it is straighforward to see that solving the vector BIE $\boldsymbol{\mathcal{B}}\boldsymbol{\phi}=\mathbf{g}$ is equivalent to solving the scalar BIE
  \begin{equation}\label{eq:vec2}
    \phi_v-\mathbf{v}\cdot \mathcal{S}(\eta)[\mathbf{u}\ \phi_v]=\mathbf{g}\cdot\mathbf{v},\quad \phi_v:=\boldsymbol{\phi}\cdot\mathbf{v}\in H^1(\Gamma).
  \end{equation}

Similar to~section \ref{sec1}, due to the invertibility of
$\mathcal{D}$ and low-rank of $\mathcal{C}$, the solution to the full
system~\eqref{eq:system_m_f} can be found by solving decoupled scalars problems.
In particular, after solving~\eqref{eq:vec2} for $\phi_v$, the Neumann
derivative of the (scattered) electric field can be recovered using
\begin{equation}\label{eq:vec2-comp}
  \partial_n E^s_j = \mathbf{e}_j \cdot \boldsymbol{\phi} = \mathbf{e}_j \cdot \left( \mathbf{g} + \mathcal{S}(\eta)[\mathbf{u}\ \phi_v] \right).
\end{equation}

\section{Spherical geometries}\label{sph}
We study in what follows the uniques solvability of the BIEs~\eqref{eq:vec1} and~\eqref{eq:vec2} in the case of spherical geometries. Because of scaling arguments, it suffices to consider the case $\Gamma=\mathbb{S}^2$. We will show in what follows that in the case when $\Gamma=\mathbb{S}^2$, the boundary integral operators featuring in equations~\eqref{eq:vec1} and~\eqref{eq:vec2} are diagonalizable in the orthormal basis of spherical harmonics, and we will compute explicitly their eigenvalues. We begin by establishing several important identities. We mention that all the identities regarding associated Legendre functions can be found in the monograph~\cite{Abramowitz}.

We begin with the following result which is a version of the result recounted in Lemma 2.4.8 in reference~\cite{Nedelec}:
\begin{lemma}\label{recursion}
  The following relations hold true for all $1\leq \ell$ and for all $|m|\leq \ell-1$
  \begin{eqnarray}\label{eq:id1}
    x_1Y_\ell^m(\hat{x})&=&\frac{1}{4}\left[-\frac{\sqrt{(\ell+m)(\ell+m-1)}}{\sqrt{\ell^2-\frac{1}{4}}}Y_{\ell-1}^{m-1}(\hat{x})+\frac{\sqrt{(\ell-m)(\ell-m-1)}}{\sqrt{\ell^2-\frac{1}{4}}}Y_{\ell-1}^{m+1}(\hat{x})\right]\nonumber\\
    &+&\frac{1}{4}\left[\frac{\sqrt{(\ell-m+1)(\ell-m+2)}}{\sqrt{(\ell+1)^2-\frac{1}{4}}}Y_{\ell+1}^{m-1}(\hat{x})-\frac{\sqrt{(\ell+m+1)(\ell+m+2)}}{\sqrt{(\ell+1)^2-\frac{1}{4}}}Y_{\ell+1}^{m+1}(\hat{x})\right]\nonumber\\
    \end{eqnarray}

    \begin{eqnarray}\label{eq:id2}
    x_2Y_\ell^m(\hat{x})&=&-\frac{i}{4}\left[\frac{\sqrt{(\ell+m)(\ell+m-1)}}{\sqrt{\ell^2-\frac{1}{4}}}Y_{\ell-1}^{m-1}(\hat{x})+\frac{\sqrt{(\ell-m)(\ell-m-1)}}{\sqrt{\ell^2-\frac{1}{4}}}Y_{\ell-1}^{m+1}(\hat{x})\right]\nonumber\\
    &+&\frac{i}{4}\left[\frac{\sqrt{(\ell-m+1)(\ell-m+2)}}{\sqrt{(\ell+1)^2-\frac{1}{4}}}Y_{\ell+1}^{m-1}(\hat{x})-\frac{\sqrt{(\ell+m+1)(\ell+m+2)}}{\sqrt{(\ell+1)^2-\frac{1}{4}}}Y_{\ell+1}^{m+1}(\hat{x})\right]\nonumber\\
    \end{eqnarray}

    and
\begin{eqnarray}\label{eq:id3}
    x_3Y_\ell^m(\hat{x})&=&\frac{1}{2}\left[\frac{\sqrt{\ell^2-m^2}}{\sqrt{\ell^2-\frac{1}{4}}}Y_{\ell-1}^m(\hat{x})+\frac{\sqrt{(\ell+1)^2-m^2}}{\sqrt{(\ell+1)^2-\frac{1}{4}}}Y_{\ell+1}^m(\hat{x})\right]
\end{eqnarray}
where $\mathbf{x}=[x_1\ x_2\ x_3]^\top$. In the case when $m=\pm \ell$, the relations above still hold with the caveat that the spherical harmonics terms $Y_L^M$ in equations~\eqref{eq:id1},~\eqref{eq:id2}, and~\eqref{eq:id3} above for which the relation $|M|\leq L$ is violated must be dropped. More precisely, we have
\begin{equation}\label{eq:id1l}
  x_1Y_\ell^\ell(\hat{x})=-\frac{1}{2}\frac{\sqrt{2\ell}}{\sqrt{2\ell+1}}Y_{\ell-1}^{\ell-1}(\hat{x})+\frac{1}{2}\frac{\sqrt{2}}{\sqrt{(2\ell+1)(2\ell+3)}}Y_{\ell+1}^{\ell-1}(\hat{x})-\frac{1}{2}\frac{\sqrt{2\ell+2}}{\sqrt{2\ell+3}}Y_{\ell+1}^{\ell+1}(\hat{x})
\end{equation}
\begin{equation}\label{eq:id2l}
  x_2Y_\ell^\ell(\hat{x})=-\frac{i}{2}\frac{\sqrt{2\ell}}{\sqrt{2\ell+1}}Y_{\ell-1}^{\ell-1}(\hat{x})+\frac{i}{2}\frac{\sqrt{2}}{\sqrt{(2\ell+1)(2\ell+3)}}Y_{\ell+1}^{\ell-1}(\hat{x})-\frac{i}{2}\frac{\sqrt{2\ell+2}}{\sqrt{2\ell+3}}Y_{\ell+1}^{\ell+1}(\hat{x})
  \end{equation}
\begin{equation}\label{eq:id3l}
  x_3Y_\ell^\ell(\hat{x})=\frac{1}{\sqrt{2\ell+3}}Y_{\ell+1}^\ell(\hat{x})
\end{equation}
with similar formulas in the case $m=-\ell$.
    
\end{lemma}
\begin{proof}
  We start with the following identity
  \[
  \frac{1}{4}\frac{\sqrt{(\ell-m)(\ell-m-1)}}{\sqrt{\ell^2-\frac{1}{4}}}Y_{\ell-1}^{m+1}(\hat{x})=\frac{1}{2(2\ell+1)}\gamma_\ell^me^{i(m+1)\varphi}\ \mathbb{P}_{\ell-1}^{m+1}(\cos{\theta})
  \]
  and
  \[
  \frac{1}{4}\frac{\sqrt{(\ell+m+1)(\ell+m+2)}}{\sqrt{(\ell+1)^2-\frac{1}{4}}}Y_{\ell+1}^{m+1}(\hat{x})=\frac{1}{2(2\ell+1)}\gamma_\ell^me^{i(m+1)\varphi}\ \mathbb{P}_{\ell+1}^{m+1}(\cos{\theta}).
  \]
  Consequently, we obtain
  \begin{eqnarray*}
    \frac{1}{4}\frac{\sqrt{(\ell-m)(\ell-m-1)}}{\sqrt{\ell^2-\frac{1}{4}}}Y_{\ell-1}^{m+1}(\hat{x})&-&\frac{1}{4}\frac{\sqrt{(\ell+m+1)(\ell+m+2)}}{\sqrt{(\ell+1)^2-\frac{1}{4}}}Y_{\ell+1}^{m+1}(\hat{x})\\
    &=&\frac{1}{2}\gamma_\ell^m e^{i(m+1)\varphi}\frac{1}{2\ell+1}\left[\mathbb{P}_{\ell-1}^{m+1}(\cos{\theta})-\mathbb{P}_{\ell+1}^{m+1}(\cos{\theta})\right].
  \end{eqnarray*}
  Using the identity
  \begin{equation}\label{eq:pml1}
    \sqrt{1-x^2}\ \mathbb{P}_\ell^m(x)=\frac{1}{2\ell+1}\left[\mathbb{P}_{\ell-1}^{m+1}(x)-\mathbb{P}_{\ell+1}^{m+1}(x)\right]
  \end{equation}
  we establish
  \begin{equation}\label{id_half1}
    \frac{1}{4}\frac{\sqrt{(\ell-m)(\ell-m-1)}}{\sqrt{\ell^2-\frac{1}{4}}}Y_{\ell-1}^{m+1}(\hat{x})-\frac{1}{4}\frac{\sqrt{(\ell+m+1)(\ell+m+2)}}{\sqrt{(\ell+1)^2-\frac{1}{4}}}Y_{\ell+1}^{m+1}(\hat{x})=\frac{1}{2}e^{i\varphi}\sin{\theta}\ Y_\ell^m(\hat{x}).
    \end{equation}
Similarly, we have
\[
  \frac{1}{4}\frac{\sqrt{(\ell+m)(\ell+m-1)}}{\sqrt{\ell^2-\frac{1}{4}}}Y_{\ell-1}^{m-1}(\hat{x})=\frac{1}{2(2\ell+1)}\gamma_\ell^me^{i(m-1)\varphi}\ (\ell+m)(\ell+m-1)\mathbb{P}_{\ell-1}^{m-1}(\cos{\theta})
  \]
  and
  \[
  \frac{1}{4}\frac{\sqrt{(\ell-m+1)(\ell-m+2)}}{\sqrt{(\ell+1)^2-\frac{1}{4}}}Y_{\ell+1}^{m-1}(\hat{x})=\frac{1}{2(2\ell+1)}\gamma_\ell^me^{i(m-1)\varphi}\ (\ell-m+1)(\ell-m+2)\mathbb{P}_{\ell+1}^{m-1}(\cos{\theta}).
  \]
  Consequently, we obtain
  \begin{eqnarray*}
    &&\frac{1}{4}\frac{\sqrt{(\ell-m+1)(\ell-m+2)}}{\sqrt{(\ell+1)^2-\frac{1}{4}}}Y_{\ell+1}^{m-1}(\hat{x})-\frac{1}{4}\frac{\sqrt{(\ell+m)(\ell+m-1)}}{\sqrt{\ell^2-\frac{1}{4}}}Y_{\ell-1}^{m-1}(\hat{x})\\
    &=&\frac{1}{2}\gamma_\ell^m e^{i(m-1)\varphi}\frac{1}{2\ell+1}\left[(\ell-m+1)(\ell-m+2)\mathbb{P}_{\ell+1}^{m-1}(\cos{\theta})-(\ell+m)(\ell+m-1)\mathbb{P}_{\ell-1}^{m-1}(\cos{\theta})\right].
  \end{eqnarray*}
 Using the identity
  \begin{equation}\label{eq:pml2}
    \sqrt{1-x^2}\ \mathbb{P}_\ell^m(x)=\frac{1}{2\ell+1}\left[(\ell-m+1)(\ell-m+2)\mathbb{P}_{\ell+1}^{m-1}(x)-(\ell+m)(\ell+m-1)\mathbb{P}_{\ell-1}^{m-1}(x)\right]
  \end{equation}
  we derive
  \begin{equation}\label{id_half2}
    \frac{1}{4}\frac{\sqrt{(\ell-m+1)(\ell-m+2)}}{\sqrt{(\ell+1)^2-\frac{1}{4}}}Y_{\ell+1}^{m-1}(\hat{x})-\frac{1}{4}\frac{\sqrt{(\ell+m)(\ell+m-1)}}{\sqrt{\ell^2-\frac{1}{4}}}Y_{\ell-1}^{m-1}(\hat{x})=\frac{1}{2}e^{-i\varphi}\sin{\theta}\ Y_\ell^m(\hat{x}).
  \end{equation}
  Clearly, combining~\eqref{id_half1} and~\eqref{id_half2}, we obtain the identities~\eqref{eq:id1} and~\eqref{eq:id2} respectively. Finally, we also have
  \[
  \frac{1}{2}\frac{\sqrt{\ell^2-m^2}}{\sqrt{\ell^2-\frac{1}{4}}}Y_{\ell-1}^m(\hat{x})=\frac{\ell+m}{2\ell+1}\gamma_\ell^me^{im\varphi}\mathbb{P}_{\ell-1}^m(\cos{\theta})
  \]
  and
  \[
  \frac{1}{2}\frac{\sqrt{(\ell+1)^2-m^2}}{\sqrt{(\ell+1)^2-\frac{1}{4}}}Y_{\ell+1}^m(\hat{x})=\frac{\ell+1-m}{2\ell+1}\gamma_\ell^me^{im\varphi}\mathbb{P}_{\ell+1}^m(\cos{\theta})
  \]
  Using the last two relations above together with the identity
  \begin{equation}\label{eq:pml3}
    (2\ell+1)x\mathbb{P}_\ell^m(x)=(\ell+m)\mathbb{P}_{\ell-1}^m(x)+(\ell+1-m)\mathbb{P}_{\ell+1}^m(x),
  \end{equation}
  the relation~\eqref{eq:id3} now follows.
  
  In the case when $m=\ell$ we have
  \[
  x_1Y_\ell^\ell(\hat{x})=\frac{1}{2}(e^{i\varphi}+e^{-i\varphi})\gamma_\ell^\ell\ e^{i\ell\varphi}\ \mathbb{P}_{\ell}^\ell(\cos{\theta})\sin{\theta}.
  \]
  We use the identity
  \begin{equation}\label{eq:plla}
    -\frac{1}{2\ell+1}\mathbb{P}_{\ell+1}^{\ell+1}(\cos{\theta})=\sin{\theta}\ \mathbb{P}_{\ell}^{\ell}(\cos{\theta})
  \end{equation}
  we derive
  \begin{eqnarray*}
    x_1Y_\ell^\ell(\hat{x})&=&-\frac{1}{2}\ \frac{1}{2\ell+1}e^{i(\ell+1)\varphi}\gamma_\ell^\ell\ \mathbb{P}_{\ell+1}^{\ell+1}(\cos{\theta})\\
    &-&\frac{1}{2}\ \frac{1}{2\ell+1}e^{i(\ell-1)\varphi}\gamma_\ell^\ell\ \mathbb{P}_{\ell+1}^{\ell+1}(\cos{\theta}).
  \end{eqnarray*}
  Since
  \[
  \gamma_\ell^\ell=(2\ell+2)\frac{\sqrt{2\ell+1}}{\sqrt{2\ell+3}}\gamma_{\ell+1}^{\ell+1}
  \]
  we obtain
  \begin{equation}\label{eq:interm1}
    x_1Y_\ell^\ell(\hat{x})=-\frac{1}{2}\ \frac{\sqrt{2\ell+2}}{\sqrt{2\ell+3}}Y_{\ell+1}^{\ell+1}(\hat{x})
    -\frac{1}{2}\ \frac{1}{2\ell+1}e^{i(\ell-1)\varphi}\gamma_\ell^\ell\ \mathbb{P}_{\ell+1}^{\ell+1}(\cos{\theta}).
  \end{equation}
 To proceed, we start with the identity
  \begin{equation}\label{eq:ident_lp1lp1}
    \mathbb{P}_{\ell+1}^{\ell+1}(\cos{\theta})=-\frac{2\ell}{\sin{\theta}}\mathbb{P}_{\ell}^{\ell}(\cos{\theta})-2\mathbb{P}_{\ell+1}^{\ell-1}(\cos{\theta})
  \end{equation}
which, using the identity~\eqref{eq:plla}, we are able to re-express in the equivalent form
\begin{equation}\label{eq:ident_lp1lp11}
    \mathbb{P}_{\ell+1}^{\ell+1}(\cos{\theta})=2\ell(2\ell-1)\mathbb{P}_{\ell-1}^{\ell-1}(\cos{\theta}) -2\mathbb{P}_{\ell+1}^{\ell-1}(\cos{\theta}).
\end{equation}
Inserting equation~\eqref{eq:ident_lp1lp11} into the formula~\eqref{eq:interm1} we obtain
\begin{eqnarray}\label{eq:interm11}
    x_1Y_\ell^\ell(\hat{x})&=&-\frac{1}{2}\ \frac{\sqrt{2\ell+2}}{\sqrt{2\ell+3}}Y_{\ell+1}^{\ell+1}(\hat{x})
    -\frac{1}{2}\ \frac{2\ell(2\ell-1)}{2\ell+1}e^{i(\ell-1)\varphi}(-1)^{\ell-1}\gamma_\ell^\ell\ \mathbb{P}_{\ell-1}^{\ell-1}(\cos{\theta})\nonumber\\
    &+&\frac{1}{2}\ \frac{2}{2\ell+1}e^{i(\ell-1)\varphi}(-1)^{\ell-1}\gamma_\ell^\ell\ \mathbb{P}_{\ell+1}^{\ell-1}(\cos{\theta}).
\end{eqnarray}
Taking into account the identities 
\[
\gamma_\ell^\ell=\frac{1}{2\ell-1}\frac{\sqrt{2\ell+1}}{\sqrt{2\ell}}\gamma_{\ell-1}^{\ell-1}\quad \gamma_\ell^\ell=\frac{1}{\sqrt{2}}\frac{\sqrt{2\ell+1}}{\sqrt{2\ell+3}}\gamma_{\ell+1}^{\ell-1}
\]
we see that equations~\eqref{eq:interm11} imply the result~\eqref{eq:id1l}. The identity~\eqref{eq:id2l} can be established similarly.

  Finally, using the identity
  \begin{equation}\label{eq:pll}
    \mathbb{P}_{\ell+1}^\ell(\cos{\theta})=(2\ell+1)\mathbb{P}_{\ell}^\ell(\cos{\theta})\cos{\theta}
  \end{equation}
  we get
  \[
  \cos{\theta}\ Y_\ell^\ell(\hat{x})=\gamma_\ell^\ell\ e^{i\ell\varphi}\ \mathbb{P}_{\ell}^\ell(\cos{\theta})\cos{\theta}=\frac{1}{\sqrt{2\ell+3}}Y_{\ell+1}^\ell(\hat{x}).
  \]
  Using the simple formula 
\[
Y_\ell^m(\hat{x})=(-1)^m \overline{Y_\ell^{-m}(\hat{x})}
\]
the identities for the case $m=-\ell$ can be produced easily. 
\end{proof}

Now we are in the position to prove our main result.
\begin{theorem}\label{thm2}
  Let $\mathcal{S}:H^p(\mathbb{S}^2)\to H^r(\mathbb{S}^2)$ be a continuous operator that is diagonalizable in the orthonormal basis $(Y_\ell^m(\hat{x}))_{0\leq\ell}^{-\ell\leq m\leq \ell}$ 
  \[
  \mathcal{S}Y_\ell^m=s_\ell Y_\ell^m,\quad for\ all\ 0\leq \ell,\ -\ell\leq m\leq \ell.
  \]
  Then the operator $\mathcal{S}_n:L^2(\mathbb{S}^2)\to L^2(\mathbb{S}^2)$ defined as
  \[
  \mathcal{S}_n\psi:=\hat{x}\cdot \mathcal{S}\left[\hat{x}\psi\right],\quad \psi\in L^2(\mathbb{S}^2)
  \]
  is also diagonalizable in the orthonormal basis $(Y_\ell^m(\hat{x}))_{0\leq\ell}^{-\ell\leq m\leq \ell}$ with eigenvalues
  \[
  \mathcal{S}_nY_\ell^m=s_{n,\ell} Y_\ell^m,\quad for\ all\ 0\leq \ell,\ -\ell\leq m\leq \ell.
  \]
  where
  \[
  s_{n,\ell}:=\begin{cases} s_{\ell+1},& \ell=0\\
  \frac{\ell}{2\ell+1}s_{\ell-1}+\frac{\ell+1}{2\ell+1}s_{\ell+1},& 1\leq\ell\end{cases}.
  \]
\end{theorem}
\begin{proof}
  We begin with the case $1\leq\ell$ and $|m|\leq \ell-1$. We rewrite equation~\eqref{eq:id1} in the form
  \[
  x_1Y_\ell^m(\hat{x})=a_\ell^mY_{\ell-1}^{m+1}(\hat{x})+b_\ell^mY_{\ell-1}^{m-1}(\hat{x})+c_\ell^mY_{\ell+1}^{m+1}(\hat{x})+d_\ell^mY_{\ell+1}^{m-1}(\hat{x})
  \]
  with coefficients $a_\ell^m,\ b_\ell^m,\ c_\ell^m$ and $d_\ell^m$ defined in same formula. Then, formula~\eqref{eq:id2} can be expressed in the form
  \[
  x_2Y_\ell^m(\hat{x})=-i\ a_\ell^mY_{\ell-1}^{m+1}(\hat{x})+i\ b_\ell^mY_{\ell-1}^{m-1}(\hat{x})-i\ c_\ell^mY_{\ell+1}^{m+1}(\hat{x})+i\ d_\ell^mY_{\ell+1}^{m-1}(\hat{x}).
  \]
  Similarly we rewrite formula~\eqref{eq:id3} in the form
  \[
  x_3Y_\ell^m(\hat{x})=e_\ell^mY_{\ell-1}^{m}(\hat{x})+f_\ell^mY_{\ell+1}^m(\hat{x})
  \]
  with coefficients $e_\ell^m$ and $f_\ell^m$ defined in the same formula. Clearly, we have
  \begin{eqnarray*}
    \mathcal{S}\left[x_1Y_\ell^m\right]&=&s_{\ell-1}\left[a_\ell^mY_{\ell-1}^{m+1}+b_\ell^mY_{\ell-1}^{m-1}\right]+s_{\ell+1}\left[c_\ell^mY_{\ell+1}^{m+1}+d_\ell^mY_{\ell+1}^{m-1}\right]\\
    \mathcal{S}\left[x_2Y_\ell^m\right]&=&s_{\ell-1}\left[-i\ a_\ell^mY_{\ell-1}^{m+1}+i\ b_\ell^mY_{\ell-1}^{m-1}\right]+s_{\ell+1}\left[-i\ c_\ell^mY_{\ell+1}^{m+1}+i\ d_\ell^mY_{\ell+1}^{m-1}\right]\\
     \mathcal{S}\left[x_3Y_\ell^m\right]&=&s_{\ell-1}e_\ell^mY_{\ell-1}^{m}+s_{\ell+1}f_\ell^mY_{\ell+1}^m.
  \end{eqnarray*}
  It follows immediately that
  \begin{eqnarray*}
    \left[\mathcal{S}_nY_\ell^m\right](\hat{x})&=&s_{\ell-1}\{\left[b_\ell^me^{i\varphi}Y_{\ell-1}^{m-1}(\hat{x})+a_\ell^m e^{-i\varphi}Y_{\ell-1}^{m+1}(\hat{x})\right]\sin{\theta}+e_\ell^mY_{\ell-1}^{m}(\hat{x})\cos{\theta}\}\\
    &+&s_{\ell+1}\{\left[d_\ell^me^{i\varphi}Y_{\ell+1}^{m-1}(\hat{x})+c_\ell^m e^{-i\varphi}Y_{\ell+1}^{m+1}(\hat{x})\right]\sin{\theta}+f_\ell^mY_{\ell+1}^{m}(\hat{x})\cos{\theta}\}.
  \end{eqnarray*}
  After some calculations we get
  \begin{eqnarray*}
    a_\ell^me^{-i\varphi}Y_{\ell-1}^{m+1}(\hat{x})&=&\frac{1}{2(2\ell+1)}\gamma_\ell^me^{im\varphi}\mathbb{P}_{\ell-1}^{m+1}(\cos{\theta})\\
    b_\ell^me^{i\varphi}Y_{\ell-1}^{m-1}(\hat{x})&=&-\frac{(\ell+m)(\ell+m-1)}{2(2\ell+1)}\gamma_\ell^me^{im\varphi}\mathbb{P}_{\ell-1}^{m-1}(\cos{\theta})\\
    e_\ell^mY_{\ell-1}^{m}(\hat{x})&=&\frac{\ell+m}{2\ell+1}\gamma_\ell^me^{im\varphi}\mathbb{P}_{\ell-1}^{m}(\cos{\theta}).
  \end{eqnarray*}
  Consequently, we obtain
  \begin{eqnarray}\label{eq:key1}
    &&\left[b_\ell^me^{i\varphi}Y_{\ell-1}^{m-1}(\hat{x})+a_\ell^m e^{-i\varphi}Y_{\ell-1}^{m+1}(\hat{x})\right]\sin{\theta}+e_\ell^mY_{\ell-1}^{m}(\hat{x})\cos{\theta}=\frac{1}{2(2\ell+1)}\gamma_\ell^me^{im\varphi}\nonumber\\
    &\times&\left[\sin{\theta}\left(\mathbb{P}_{\ell-1}^{m+1}(\cos{\theta})-(\ell+m)(\ell+m-1)\mathbb{P}_{\ell-1}^{m-1}(\cos{\theta})\right)+2(\ell+m)\cos{\theta}\ \mathbb{P}_{\ell-1}^m(\cos{\theta})\right].\nonumber\\
  \end{eqnarray}
  We make use of the following identities
  \begin{equation}\label{eq:pml4}
    2m\mathbb{P}_\ell^m(x)=-\sqrt{1-x^2}\left[\mathbb{P}_{\ell-1}^{m+1}(x)+(\ell+m)(\ell+m-1)\mathbb{P}_{\ell-1}^{m-1}(x)\right]
  \end{equation}
  and
  \begin{equation}\label{eq:pml5}
    \sqrt{1-x^2}\mathbb{P}_{\ell-1}^{m+1}(x)=(\ell-m)\mathbb{P}_\ell^m(x)-(\ell+m)x\mathbb{P}_{\ell-1}^m(x)
  \end{equation}
  and derive
  \begin{equation}\label{eq:key11}
    \sqrt{1-x^2}\left[\mathbb{P}_{\ell-1}^{m+1}(x)-(\ell+m)(\ell+m-1)\mathbb{P}_{\ell-1}^{m-1}(x)\right]+2(\ell+m)x\mathbb{P}_{\ell-1}^m(x)=2\ell \mathbb{P}_\ell^m(x).
  \end{equation}
  Using the identity~\eqref{eq:key11} in equation~\eqref{eq:key1}, we derive
  \begin{equation}\label{eq:key111}
    \left[b_\ell^me^{i\varphi}Y_{\ell-1}^{m-1}(\hat{x})+a_\ell^m e^{-i\varphi}Y_{\ell-1}^{m+1}(\hat{x})\right]\sin{\theta}+e_\ell^mY_{\ell-1}^{m}(\hat{x})\cos{\theta}=\frac{\ell}{2\ell+1}Y_\ell^m(\hat{x}).
  \end{equation}
  Similarly, we obtain
  \begin{eqnarray*}
    c_\ell^me^{-i\varphi}Y_{\ell+1}^{m+1}(\hat{x})&=&-\frac{1}{2(2\ell+1)}\gamma_\ell^me^{im\varphi}\mathbb{P}_{\ell+1}^{m+1}(\cos{\theta})\\
    d_\ell^me^{i\varphi}Y_{\ell+1}^{m-1}(\hat{x})&=&\frac{(\ell-m+1)(\ell-m+2)}{2(2\ell+1)}\gamma_\ell^me^{im\varphi}\mathbb{P}_{\ell+1}^{m-1}(\cos{\theta})\\
    f_\ell^mY_{\ell+1}^{m}(\hat{x})&=&\frac{\ell+1-m}{2\ell+1}\gamma_\ell^me^{im\varphi}\mathbb{P}_{\ell+1}^{m}(\cos{\theta}).
  \end{eqnarray*}
  Consequently, we obtain
  \begin{eqnarray}\label{eq:key2}
    &&\left[d_\ell^me^{i\varphi}Y_{\ell+1}^{m-1}(\hat{x})+c_\ell^m e^{-i\varphi}Y_{\ell+1}^{m+1}(\hat{x})\right]\sin{\theta}+f_\ell^mY_{\ell+1}^{m}(\hat{x})\cos{\theta}=\frac{1}{2(2\ell+1)}\gamma_\ell^me^{im\varphi}\nonumber\\
    &\times&\left[\sin{\theta}\left(-\mathbb{P}_{\ell+1}^{m+1}(\cos{\theta})+(\ell-m+1)(\ell-m+2)\mathbb{P}_{\ell+1}^{m-1}(\cos{\theta})\right)+2(\ell+1-m)\cos{\theta}\ \mathbb{P}_{\ell+1}^m(\cos{\theta})\right].\nonumber\\
  \end{eqnarray}
  We make use of the following identities
  \begin{equation}\label{eq:pml6}
    -2m\mathbb{P}_\ell^m(x)=\sqrt{1-x^2}\left[\mathbb{P}_{\ell+1}^{m+1}(x)+(\ell-m+1)(\ell-m+2)\mathbb{P}_{\ell+1}^{m-1}(x)\right]
  \end{equation}
  and
  \begin{equation}\label{eq:pml7}
    \sqrt{1-x^2}\mathbb{P}_{\ell+1}^{m+1}(x)=(\ell+1-m)x\mathbb{P}_{\ell+1}^m(x)-(\ell+1+m)\mathbb{P}_\ell^m(x)
  \end{equation}
  and derive
  \begin{equation}\label{eq:key22}
    \sqrt{1-x^2}\left[-\mathbb{P}_{\ell+1}^{m+1}(x)+(\ell-m+1)(\ell-m+2)\mathbb{P}_{\ell+1}^{m-1}(x)\right]+2(\ell+1-m)x\mathbb{P}_{\ell+1}^m(x)=2(\ell+1) \mathbb{P}_\ell^m(x).
  \end{equation}
  Using the identity~\eqref{eq:key22} in equation~\eqref{eq:key2}, we derive
  \begin{equation}\label{eq:key222}
    \left[d_\ell^me^{i\varphi}Y_{\ell+1}^{m-1}(\hat{x})+c_\ell^m e^{-i\varphi}Y_{\ell+1}^{m+1}(\hat{x})\right]\sin{\theta}+f_\ell^mY_{\ell+1}^{m}(\hat{x})\cos{\theta}=\frac{\ell+1}{2\ell+1}Y_\ell^m(\hat{x}).
  \end{equation}
  Finally, combining~\eqref{eq:key111} and~\eqref{eq:key222} we derive
  \begin{equation}\label{eq:eigenvalue1}
    \mathcal{S}_nY_\ell^m=\left[\frac{\ell}{2\ell+1}s_{\ell-1}+\frac{\ell+1}{2\ell+1}s_{\ell+1}\right]Y_\ell^m,\quad |m|\leq\ell-1.
  \end{equation}

  In the case when $1\leq \ell=m$, we rewrite equations~\eqref{eq:id1l}--\eqref{eq:id3l} in the form
  \[
  x_1Y_\ell^\ell(\hat{x})=b_\ell Y_{\ell-1}^{\ell-1}(\hat{x})+c_\ell Y_{\ell+1}^{\ell+1}(\hat{x}) +d_\ell Y_{\ell+1}^{\ell-1}(\hat{x})
  \]
  and
  \[
  x_2Y_\ell^\ell(\hat{x})=ib_\ell Y_{\ell-1}^{\ell-1}(\hat{x})-ic_\ell Y_{\ell+1}^{\ell+1}(\hat{x}) +id_\ell Y_{\ell+1}^{\ell-1}(\hat{x})
  \]
  while
  \[
  x_3Y_\ell^\ell(\hat{x})=e_\ell Y_{\ell+1}^{\ell}(\hat{x}).
  \]
  A simple calculation leads to the following relation
  \begin{eqnarray*}
    [S_nY_\ell^\ell](\hat{x})&=&\left(\sin{\theta}(d_\ell e^{i\varphi}Y_{\ell+1}^{\ell-1}(\hat{x}) + c_\ell e^{-i\varphi}Y_{\ell+1}^{\ell+1}(\hat{x}))+e_\ell Y_{\ell+1}^\ell(\hat{x})\cos{\theta}\right)s_{\ell+1}\\
    &+& (b_\ell e^{i\varphi}Y_{\ell-1}^{\ell-1}(\hat{x})\sin{\theta})s_{\ell-1}.
  \end{eqnarray*}
  It is straightforward to derive the following identities
  \begin{eqnarray*}
    c_\ell e^{-i\varphi}Y_{\ell+1}^{\ell+1}(\hat{x})&=&-\frac{1}{2(2\ell+1)} e^{i\ell\varphi}\gamma_\ell^\ell \mathbb{P}_{\ell+1}^{\ell +1}(\cos{\theta})\\
    d_\ell e^{i\varphi}Y_{\ell+1}^{\ell-1}(\hat{x})&=&\frac{1}{2\ell+1} e^{i\ell\varphi}\gamma_\ell^\ell \mathbb{P}_{\ell+1}^{\ell-1}(\cos{\theta})\\
    e_\ell Y_{\ell+1}^\ell(\hat{x})\cos{\theta}&=&\frac{1}{2\ell+1} e^{i\ell\varphi}\gamma_\ell^\ell \mathbb{P}_{\ell+1}^{\ell}(\cos{\theta})\ \cos{\theta}.
  \end{eqnarray*}
  Therefore, we obtain
  \begin{eqnarray*}
    \sin{\theta}&\times&(d_\ell e^{i\varphi}Y_{\ell+1}^{\ell-1}(\hat{x}) + c_\ell e^{-i\varphi}Y_{\ell+1}^{\ell+1}(\hat{x}))+e_\ell Y_{\ell+1}^\ell(\hat{x})\cos{\theta}=\frac{1}{2\ell+1} e^{i\ell\varphi}\gamma_\ell^\ell\\
    &\times&\left(\mathbb{P}_{\ell+1}^{\ell}(\cos{\theta})\cos{\theta}-\frac{1}{2}\mathbb{P}_{\ell+1}^{\ell+1}(\cos{\theta})\sin{\theta}+\mathbb{P}_{\ell+1}^{\ell-1}(\cos{\theta})\sin{\theta}\right)
  \end{eqnarray*}
  Using identity~\eqref{eq:ident_lp1lp1} we obtain
  \[
  \mathbb{P}_{\ell+1}^{\ell}(\cos{\theta})\cos{\theta}-\frac{1}{2}\mathbb{P}_{\ell+1}^{\ell+1}(\cos{\theta})\sin{\theta}+\mathbb{P}_{\ell+1}^{\ell-1}(\cos{\theta})\sin{\theta}=\mathbb{P}_{\ell+1}^{\ell}(\cos{\theta})\cos{\theta}-\mathbb{P}_{\ell+1}^{\ell+1}(\cos{\theta})\sin{\theta}-\ell\mathbb{P}_\ell^\ell(\cos{\theta}).
  \]
  Now using the identity~\eqref{eq:plla} together with the identity~\eqref{eq:pll}
  we get
  \[
  \mathbb{P}_{\ell+1}^{\ell}(\cos{\theta})\cos{\theta}-\frac{1}{2}\mathbb{P}_{\ell+1}^{\ell+1}(\cos{\theta})\sin{\theta}+\mathbb{P}_{\ell+1}^{\ell-1}(\cos{\theta})\sin{\theta}=(\ell+1)\mathbb{P}_\ell^\ell(\cos{\theta})
  \]
  and thus
  \begin{equation}\label{eq:final_id_1}
    \sin{\theta}(d_\ell e^{i\varphi}Y_{\ell+1}^{\ell-1}(\hat{x}) + c_\ell e^{-i\varphi}Y_{\ell+1}^{\ell+1}(\hat{x}))+e_\ell Y_{\ell+1}^\ell(\hat{x})\cos{\theta}=\frac{\ell+1}{2\ell+1}Y_\ell^\ell(\hat{x}).
  \end{equation}
  Finally, using the identity~\eqref{eq:plla} we obtain
  \begin{equation}\label{eq:final_id_2}
    b_\ell e^{i\varphi}Y_{\ell-1}^{\ell-1}(\hat{x})\sin{\theta}=\frac{\ell}{2\ell+1}Y_\ell^\ell(\hat{x}).
  \end{equation}
  Combining~\eqref{eq:final_id_1} and~\eqref{eq:final_id_2} we get
  \begin{equation}\label{eq:eigenvalue1}
    \mathcal{S}_nY_\ell^\ell=\left[\frac{\ell}{2\ell+1}s_{\ell-1}+\frac{\ell+1}{2\ell+1}s_{\ell+1}\right]Y_\ell^\ell,\quad 1\leq\ell.
  \end{equation}
Similar calculations lead to the same eignevalue formula in the case $m=-\ell$ and $\ell>0$. Finally, in the case $\ell=0$, we use the identities
  \begin{eqnarray}\label{eq:spherical_harm}
    x_1&=&\sqrt{\frac{2\pi}{3}}\left(Y_1^{-1}(\theta,\varphi)-Y_1^1(\theta,\varphi)\right)\nonumber\\
    x_2&=&i\sqrt{\frac{2\pi}{3}}\left(Y_1^{-1}(\theta,\varphi)+Y_1^1(\theta,\varphi)\right)\nonumber\\
    x_3&=&\sqrt{\frac{2\pi}{3}}Y_1^{0}(\theta,\varphi).
  \end{eqnarray}
  to derive
  \begin{equation}\label{eq:eigenvalue2}
    \mathcal{S}_nY_0^0=s_{1}Y_0^0.
    \end{equation}
\end{proof}

Given the result in Theorem~\ref{thm2}, the spectral properties of the boundary integral operators that enter formulations~\eqref{eq:vec1} and~\eqref{eq:vec2} in the case when $\Gamma=\mathbb{S}^2$ hinge on the spectral properties of the boundary integral operator $\mathcal{S}(\eta)=i\eta\mathcal{D}^{-1}(\eta)\mathcal{N}(\eta^{-1})$. Given that all of the four boundary integral operators associated with the Helmholtz Calder\'on calculus are diagonalizable in the case of the unit sphere, it turns out that so is the operator $\mathcal{S}(\eta)$. Indeed, given the spectral relations~\cite{turc3} which hold for all $0\leq \ell$ and $-\ell\leq m\leq\ell$
\begin{eqnarray*}
  SY_\ell^m&=&ikj_\ell(k)h_\ell^{(1)}(k)Y_\ell^m\\
  KY_\ell^m&=&\left(-\frac{1}{2}+ik^2j'_\ell(k)h_\ell^{(1)}(k)\right)Y_\ell^m=K^\top Y_\ell^m\\
  NY_\ell^m&=&ik^3j'_\ell(k)(h_\ell^{(1)})'(k)Y_\ell^m\\
\end{eqnarray*}
a simple calculations reveals
\begin{equation}\label{eq:sell}
 \mathcal{S}(\eta)Y_\ell^m=s_\ell(k) Y_\ell^m\qquad  s_\ell(k)=k\frac{(h_\ell^{(1)}(k))'}{h_\ell^{(1)}(k)}
\end{equation}
where we make the convention that $s_{-1}(k):=0$. Remarkably, the eigenvalues of the operator $\mathcal{S}(\eta)$ do not depend on the coupling parameter $\eta$. The role of the coupling parameter is to ensure that the operator $\mathcal{S}(\eta)$ is properly defined for all wavenumbers $k>0$. Using the result in Theorem~\ref{thm2}, we obtain the following spectral properties for the scalar Woodbury boundary integral operators that enter formulations~\eqref{eq:vec1} and~\eqref{eq:vec2} in the case of a unit sphere scatterer:
\begin{equation}\label{eq:eig_11}
  Y_\ell^m+\mathbf{n}\cdot\mathcal{S}(\eta)\left[\frac{1}{2}\mathcal{H}^{-1}\mathbf{n}\ Y_\ell^m\right]=\left(1+\frac{\ell}{2(2\ell+1)}s_{\ell-1}+\frac{\ell+1}{2(2\ell+1)}s_{\ell+1}\right)Y_\ell^m
\end{equation}
and
\begin{equation}\label{eq:eig_22}
  Y_\ell^m-\begin{bmatrix}-\mathbf{n} \\ 1\end{bmatrix}\cdot\mathcal{S}(\eta)\left[\begin{bmatrix}\mathbf{n}\\ 1\end{bmatrix}Y_\ell^m\right]=\left(1+\frac{\ell}{(2\ell+1)}s_{\ell-1}+\frac{\ell+1}{(2\ell+1)}s_{\ell+1}-s_\ell\right)Y_\ell^m
\end{equation}

We are now in the position to prove our main result
\begin{theorem}\label{thm1} Both boundary integral equation formulations~\eqref{eq:vec1} and~\eqref{eq:vec2} are uniquely solvable in the case when $\Gamma=\mathbb{S}^2$ for all wavenumbers $k>0$.
\end{theorem}
\begin{proof} According to equations~\eqref{eq:eig_11}, the unique solvability of the equation~\eqref{eq:vec1} in the case of the unit sphere is equivalent to showing that
\[
2+\frac{\ell}{2\ell+1}s_{\ell-1}+\frac{\ell+1}{2\ell+1}s_{\ell+1}\neq 0,\ 0\leq \ell
\]
where $s_\ell$ were defined inequations~\eqref{eq:sell}. According to reference~\cite{Nedelec} (formulas (2.6.23) and (2.6.24)) we have that
\[
s_\ell(k)=-\frac{p_\ell(k)}{q_\ell(k)}+i\frac{k}{q_\ell(k)}
\]
where
\[
q_\ell(k)=1+\alpha_1^\ell\frac{1}{k^2}+\ldots+\alpha_\ell^\ell\frac{1}{k^{2\ell}}
\]
and
\[
p_\ell(k)=1+2\alpha_1^\ell\frac{1}{k^2}+\ldots+(\ell+1)\alpha_\ell^\ell\frac{1}{k^{2\ell}}
\]
with
\[
\alpha_m^\ell=\beta_m^\ell \beta_m^m,\qquad \beta_m^\ell:=\frac{(m+\ell)!}{m!(\ell-m)!2^m}.
\]
It can be clearly seen that given that $\alpha_m^\ell>0$, then
\begin{equation}\label{eq:ned_decisive}
-\Re(s_\ell(k))>1,\qquad \Im(s_\ell(k))> 0\qquad \ell\geq 0\qquad k>0.
\end{equation}
If we assume by contradiction that
\[
2=\frac{\ell}{2\ell+1}(-s_{\ell-1})+\frac{\ell+1}{2\ell+1}(-s_{\ell+1})
\]
for $\ell\geq 1$, then we get
\[
2=\frac{\ell}{2\ell+1}\Re(-s_{\ell-1})+\frac{\ell+1}{2\ell+1}\Re(-s_{\ell+1})
\]
and
\[
0=\frac{\ell}{2\ell+1}\Im(s_{\ell-1})+\frac{\ell+1}{2\ell+1}\Im(s_{\ell+1})
\]
which would contradict inequalities~\eqref{eq:ned_decisive}. In the case $\ell=0$, the eigenvalue corresponding to zero mode equals $1+s_1/2$ whose imaginary part is again nozero according to inequalities~\eqref{eq:ned_decisive}.

Given equations~\eqref{eq:eig_11}, the unique solvability of the formulation~\eqref{eq:vec2}  is equivalent to proving 
\[
1+\frac{\ell}{2\ell+1}s_{\ell-1}+\frac{\ell+1}{2\ell+1}s_{\ell+1}-s_\ell\neq 0,\ 0\leq \ell.
\]
If we assume by contradiction that
\begin{equation}\label{eq:contra}
1+\frac{\ell}{2\ell+1}s_{\ell-1}+\frac{\ell+1}{2\ell+1}s_{\ell+1}=s_\ell,
\end{equation}
and use the recursion relations (equation (2.6.18) in~\cite{Nedelec})
\[
(s_{\ell-1}-(\ell-1))(s_\ell+(\ell+1))=-k^2
\]
and
\[
(s_\ell-\ell)(s_{\ell+1}+(\ell+2))=-k^2
\]
we obtain that the relation~\eqref{eq:contra} is equivalent to
\begin{equation}\label{eq:contra1}
-\frac{k^2}{2\ell+1}\left(\frac{\ell}{s_\ell+\ell+1}+\frac{\ell+1}{s_\ell-\ell}\right)=s_\ell+1.
\end{equation}
We note that given that $\Im{s_\ell}>0$, the denominators of the quantities featured in equation~\eqref{eq:contra1} do not vanish. Furthermore, equation~\eqref{eq:contra1} can be expressed in the equivalent form
\begin{equation}\label{eq:contra2}
-\frac{k^2(s_\ell+1)}{(s_\ell+\ell+1)(s_\ell-\ell)}=s_\ell+1.
\end{equation}
Now, we have that
\[
s_\ell=-\frac{p_\ell}{q_\ell}+i\ \frac{k}{q_\ell}
\]
and thus $s_\ell+1\neq 0$. Thus, equation~\eqref{eq:contra2} is equivalent to
\begin{equation}\label{eq:contra3}
  (s_\ell+\ell+1)(s_\ell-\ell)=-k^2.
  \end{equation}
Taking the imaginary part of equation~\eqref{eq:contra3} we obtain that
\[
\frac{p_\ell}{q_\ell}=\frac{1}{2},
\]
which is clearly a contradiction given that we have already established the fact that
\[
\frac{p_\ell}{q_\ell}>1.
\]
In conclusion, relation~\eqref{eq:contra} cannot hold.
\end{proof}

Having established that neither of the eigenvalues in equations~\eqref{eq:eig_11} and ~\eqref{eq:eig_22} vanish for all values of the wavenumber $k>0$, we investigate in what follows their asymptotic properties for fixed $k$ and large values of the index $\ell$. To this end, using the definition of the spherical Hankel functions and known recurrence formulas for the derivatives of Hankel functions, we express the quantities $s_\ell$ in the following equivalent form
\begin{equation}\label{eq:diff_form}
  s_\ell=\ell-ka_\ell,\qquad a_\ell:=\frac{H_{\ell+\frac{3}{2}}^{(1)}(k)}{H_{\ell+\frac{1}{2}}^{(1)}(k)},
\end{equation}
where $H_\ell^{(1)}$ denote the Hankel functions of the first kind. Using the asymptotic formula (10.19.2) in~\cite{Abramowitz}, valid for fixed argument $z$ and index $\nu\to\infty$
\[
H_\nu^{(1)}(z)= -i\sqrt{\frac{2}{\pi\nu}}\left(\frac{ez}{2\nu}\right)^{-\nu}[1+\mathcal{O}(\nu^{-1})],
\]
we derive immediately the following asymptotic relation valid for a fixed wavenumber $k$
\begin{equation}\label{eq:asympt_al}
  a_\ell=\frac{2\ell+a}{k}+\mathcal{O}(\ell^{-1}),\qquad s_\ell=-(\ell+a)+\mathcal{O}(\ell^{-1}),\qquad \ell\to\infty.
\end{equation}
First, a direct application of formula~\eqref{eq:asympt_al} leads to the following asymptotic relation of the eigenvalues in equation~\eqref{eq:eig_11}
\begin{equation}\label{eq:asympt1}
  1+\frac{\ell}{2(2\ell+1)}s_{\ell-1}+\frac{\ell+1}{2(2\ell+1)}s_{\ell+1}\sim-\frac{2\ell^2+3\ell+2}{2(2\ell+1)}\sim-\frac{\ell}{2},\quad \ell\to\infty.
\end{equation}
Therefore, we conclude from the asymptotic formulas~\eqref{eq:asympt1} that the BIE formulation that involves the mean curvature in equations~\eqref{eq:eig_11} is of the first kind in the case of spherical geometries since their spectra accumulate at infinity. The situation is very different in the case of the BIE formulation featured in equations~\eqref{eq:eig_22}. Indeed, using the asymptotic relations~\eqref{eq:asympt_al} we obtain
\begin{equation}\label{eq:asympt2}
1+\frac{\ell}{(2\ell+1)}s_{\ell-1}+\frac{\ell+1}{(2\ell+1)}s_{\ell+1}-s_\ell=1+\mathcal{O}(\ell^{-1}),\ \ell\to\infty
\end{equation}
which shows that eigenvalues in the formulation~\eqref{eq:eig_22} accumulate at 1 for fixed $k$ and large values of the index $\ell$, and therefore the scalar Woodbury BIE associated with the augmented field-only BIE is of the second kind, at least for spherical geometries.

\section{Conclusions}

We presented an analysis of the spectral properties of combined field-only boundary integral equation formulations of Maxwell scattering problems for PEC spheres. The proof of well-posedness of these formulations as well as their numerical implementation using Nystr\"om discretizations are currently underway. 

\section*{Acknowledgments}
Catalin Turc gratefully acknowledge support from NSF through contracts DMS-1908602.
\bibliography{biblioLayer}

\begin{thebibliography}{80}

 \bibitem{Abramowitz} Abramowitz, M., and I.Stegun, ``Handbook of mathematical functions with formulas, graphs, and mathematical tables'', vol. 55 of National Bureau of Standards Applied Mathematics Series. For sale by the Superintendent of Documents, U.S. Government Printing Office, Washington D.C., 1964.
 

\bibitem {AlougesLevadoux} Alouges, F., Borel, S., and D. Levadoux,
   ``A stable well-conditioned integral equation for electromagnetism
   scattering", in \emph{J. Comput. Appl. Math., 204}, 2007,
   pp. 440--451.

   \bibitem{	HarringtonMautz} Harrington, R., Mautz J., \emph{H-field, E-field and combined field solution for conducting bodies of revolution}, Arch. Elek. Uber. (AEU) 32 (1978), 4, 157-164
 
 \bibitem{Yuffa1} Yuffa, A. J., and J. Markkanen, \emph{A 3-D tensorial integral formulation of scattering containing intriguing relations}, IEEE. Trans. Antennas Propag., vol 66, no 10, pp 5274-5281

 \bibitem{Yuffa2} D. Y. C. Chan, A. J. Yuffa, E. Klaseboer, and Q. Sun, \emph{Efficient field-only surface integral equations for electromagnetics}, JOSA A, vol 37, 2, 2020, 276--283.  

\bibitem{Klaseboer} Klaseboer, E., Q. SUn, and D. Y. Chan, \emph{Non-singular field-only surface integral equations for electromagnetic scattering}, IEEE. Trans. Antennas Propag., vol 65, no 2, pp 972-977 

 \bibitem {BrackhageWerner} Brackhage, H., and P. Werner, \emph{Uber das Dirichletsche
 Aussenraumproblem fur die Helmholtsche Schwingungsgleichung},
 Arch.Math, 16, 325-329, 1965.
 
\bibitem{br-turc} Bruno, O., Elling, T., Turc, C., \emph{Regularized integral equations and fast high-order solvers for sound-hard acoustic scattering problems}, Int. J. Eng. Math., 2011.
 
\bibitem {BurtonMiller} Burton, A. J., and G. F. Miller, ``The
   application of integral equation methods to the numerical solution
   of some exterior boundary-value problems'', in \emph{Proc. Royal
     Soc. London}, vol 323, 1971, pp. 201--210.
 
\bibitem{Burton} Burton, A., Numerical solution of acoustic radiation problems, \emph{NPL Contract Rept. OC5/S35 National Physical Laboratory, Teddington, Middlesex}, (1976).
 
 
\bibitem{turc3} Boubendir, Y., Turc, C., \emph{Well-conditioned boundary integral equation formulations for the solution of high-frequency electromagnetic scattering problems}, Computer \& Mathematics with Applications, 67 (10), 2014, 1772--1805.
 

\bibitem{HsiaoKleinman} Hsiao, G. C., and R. E. Kleinman R~E,   ``Mathematical foundations for error estimation in numerical solutions of integral equations in electromagnetics'', \emph {IEEE Trans. Antennas and Propag. 45}, 1997, pp 316--328.
 

\bibitem{epstein} C. Epstein and L. Greengard, \emph{Debye sources and the numerical solution of the time harmonic Maxwell equations}, Comm. Pure Appl. Math., 63 (4), 413--463, 2010.


\bibitem{Nedelec} N\'ed\'elec, J. C., \emph{Acoustic and electromagnetic equations}, Springer Verlag, New York, 2001.


\bibitem{PoggioMiller} Poggio A. J., and E. K. Miller, ``Integral
equation solutions of three-dimensional scattering problems'', in
\emph{Computer Techniques for Electromagnetics}, R. Mittra, Ed. New York: Pergamon, 1973, ch. 4, pp. 201--210.
  

\bibitem{Woodbury} M. Woodbury, \emph{Inverting modified matrices}, Memorandum Rept. 42, Statistical Research Group, Princeton University, princeton, NJ, 1950, 4pp.




\end{thebibliography}

\end{document}